\documentclass{article}
\usepackage[english]{babel}
\usepackage[letterpaper,top=2cm,bottom=2cm,left=3cm,right=3cm,marginparwidth=1.75cm]{geometry}
\usepackage{amsmath,amssymb}
\usepackage{graphicx}
\usepackage[colorlinks=true, allcolors=cyan]{hyperref}
\usepackage{lipsum} 
\usepackage{authblk}
\usepackage{url}
\usepackage{hyperref}
\usepackage{color}
\usepackage[normalem]{ulem} 
\usepackage{pbox} 
\usepackage{makecell}
\usepackage{tabularx}
\usepackage{pbox}
\usepackage{tabulary}
\usepackage{array}
\usepackage{hyperref}
\usepackage{rotating}
\usepackage{lineno}
\usepackage{graphicx}
\usepackage{subcaption}
\usepackage{multirow} 

\renewcommand{\arraystretch}{1.5}

\providecommand{\keywords}[1]
{
	
	\textbf{{Keywords:}} #1
}

\title{Competitive tumor growth modeling and optimal radiotherapy control via logistic equations}

\author[1,2,\thanks{corresponding author, email: \href{mailto:javier.lopez.pedrares@uvigo.gal}{javier.lopez.pedrares@uvigo.gal}}]{\underline{Javier L\'opez-Pedrares}}
\author[2]{Alba López-Rivas}
\author[2]{Raquel Romero-Lorenzo}
\author[4]{Jacobo Guiu-Souto}
\author[2,3]{Alberto P. Mu\~nuzuri}

\affil[1]{Department of Mathematics, Universidade de Vigo, 36310 Vigo, Spain}
\affil[2]{Group of Nonlinear Physics, Universidade de Santiago de Compostela, 15782 Santiago de Compostela, Spain}
\affil[3]{Galician Center for Mathematical Research and Technology (CITMAga), 15782 Santiago de Compostela, Spain}
\affil[4]{Department of Medical Physics, Centro Oncolóxico de Galicia, 15009 A Coruña, Spain}

\begin{document}
	\date{}
	\maketitle
	\begin{abstract}
		
	The uncontrolled proliferation of cancer cells and their interaction with healthy tissue poses a major challenge in oncology. This manuscript develops and analyzes mathematical models that describe tumor response to radiotherapy by incorporating the Linear–Quadratic model for cell survival. To improve therapeutic efficiency, the theory of optimal control is introduced on a system of coupled differential equations, allowing for the comparison of constant versus optimized radiation strategies. The analytical study of these models provides insights into the expected dynamics under different treatment scenarios, while numerical simulations validate the theoretical results and highlight the benefits of optimal control in reducing tumor burden with minimized collateral damage.
		
	\end{abstract}
	

    \section{Introduction}

In terms of mortality, cancer is the second leading cause of death worldwide after cardiovascular diseases, accounting for approximately 9.7 million deaths in 2022 \cite{Bray2024GLOBOCAN}. However, current trends suggest it may become the leading cause of death in the coming years. In Spain, the incidence of cancer, defined as the number of newly diagnosed cases, was estimated to reach 296,103 in 2025, representing a 3.3\% increase with respect to the previous year \cite{SEOM2025}.

This clinical burden has motivated the development of mathematical models to study tumor growth dynamics in an attempt to define novel strategies to control and eventually eliminate it. Early mathematical models of tumor growth appeared in the 1970s and were mainly population-based \cite{NortonSimon1977}. However, pioneering efforts in radiobiology date back even further, such as the work of Schwartz in 1961, which analyzed tumor response to radiation through mathematical frameworks \cite{Schwartz1961}. Within this framework, the Gompertz model became a standard for describing the growth kinetics of various tumors, as it accounts for the observed decrease in growth rate as the tumor volume increases \cite{Laird1964}. During the 1990s, models aimed at describing the spatial growth of tumors began to be developed, following the pioneering work by Greenspan in this field in 1972 \cite{Greenspan1972ModelsFT}. These studies led to approaches spanning simple diffusive formulations to more sophisticated multiphase mechanical models.

In recent years, research in mathematical oncology has experienced a significant increase, driven by the integration of quantitative and computational approaches into biomedical research \cite{dOnofrioGandolfi2013}. Specifically, the Linear-Quadratic (LQ) model has been widely adopted to quantify the biological effect of radiotherapy, establishing a mathematical relationship between the radiation dose and the cell survival fraction \cite{hall2018radiobiology, fowler1989linear}. These models provide a rigorous framework for understanding free tumor growth and optimizing existing treatments, as well as exploring new therapeutic strategies \cite{Altrock2015,McMahon2018LQ,Agur2014}. Many of these works have been validated, and their reliability is well supported by empirical data \cite{Benzekry2014, Hidrovo2017}.

Although tumor dynamics are highly complex, it can generally be stated that, in the absence of treatment, tumors tend to progressively dominate the host organism, potentially leading to patient death \cite{Zolkind2021UntreatedHNC}. Nevertheless, tumor growth is not unbounded, since the human body is subject to physical and metabolic constraints. In particular, the energetic cost associated with cellular proliferation leads to limitations in oxygen and nutrient supply, resulting in the formation of hypoxic regions within the tumor core and, consequently, leading to the existence of a growth limit \cite{Schattler2015OptimalControlCancer}.

Currently, the main therapeutic modalities used in cancer treatment include chemotherapy, surgery, and radiotherapy \cite{CiudadPlatero2003}. Among these, radiotherapy represents a substantial component of hospital activity in oncology. Given the high incidence of cancer in the population, the number of patients receiving radiotherapy each year is considerable, estimated at approximately 120000 in Spain according to the Spanish Society of Radiation Oncology (SEOR), which reports that between 50\% and 60\% of oncology patients undergo this type of treatment at some point during their clinical management \cite{SEOR2021LibroBlanco}. This scale highlights the clinical relevance of improving the understanding of tumor evolution during radiotherapy and motivates the development of simulation-based approaches, including optimal control frameworks, that may support clinical decision-making within the context of personalized medicine. 

The robustness and effectiveness of optimal control techniques have been well demonstrated in other biomedical contexts, particularly in epidemiology, where they have been widely used to design vaccination , treatment, and containment strategies for infectious diseases. For instance, at the end of the $20^{th}$ century, Kirschner, Lenhart and Serbin modeled immune-response treatment strategies using optimal control methods \cite{kirschner1997}. Shortly thereafter, Behncke developed an optimal control framework for HIV therapy scheduling \cite{behncke2000}. Since then, a wide range of epidemiological problems have been addressed within this framework, including SIR type epidemic models, vaccination policy design, and, more recently, intervention strategies during the COVID-19 pandemic \cite{lenhart2007, ledzewicz2011, sharomi2017optimal, DjidjouCOVID}. These studies show that optimal control can provide realistic strategies in complex biological systems supported by empirical data.

Therefore, applying this perspective to oncology is a natural and promising step. Although optimal control has already been explored in chemotherapy scheduling and treatment optimization \cite{swan1980, martin1994, swierniak2016cell, ledzewicz2012tumor}, its integration into comprehensive tumor growth and response frameworks, particularly in radiotherapy, remains limited. In this context, optimal control theory provides a powerful framework for reinterpreting tumor growth and therapy as a dynamic optimization problem rather than a purely descriptive biological process. Instead of only predicting tumor evolution, it enables the identification of treatment protocols that minimize tumor burden while accounting for toxicity, resistance, and patient-specific constraints. Radiotherapy dose, timing, and fractionation can thus be treated as control variables acting on a nonlinear dynamical system, allowing a systematic exploration of therapeutic trade-offs and supporting the design of adaptive, personalized strategies.

This article focuses on the effect of introducing control functions into classical competition models to describe radiotherapy’s impact on cell proliferation. A key contribution of this work is the transition from simplified closed-system models to an open-system framework that reflects the metabolic struggle between healthy and malignant cells. We study how therapeutic commands can break the dominance of cancerous populations in the free competition model, effectively reversing the biological ``sink” of the untreated disease.

The manuscript is organized as follows. In Section \ref{sec:Model_No_Control}, we present different models of cancer spread, progressing from single-population growth to a coupled competition system, including an analytical study of the expected stability and numerical simulations. In Section \ref{sec:Model_Control}, we introduce the optimal control problem associated with the model, presenting both an analytical study for constant treatment and numerical simulations for dynamic dose optimization. The final section provides the discussion and conclusions, highlighting the potential for future clinical applications.

\section{Control-free tumor growth model}
\label{sec:Model_No_Control}

The mathematical description of tumor progression and its interaction with surrounding tissues is essential to predict clinical outcomes and optimize therapeutic interventions \cite{enderling2014mathematical}. In this work, the biological environment is modeled by considering two distinct cell populations: cancerous cells, $C$, and healthy cells, $H$. The temporal evolution of these populations is governed by their intrinsic growth dynamics and their specific physiological constraints.

While tumor cells are characterized by unregulated proliferation and an evasion of homeostatic signals, healthy tissues maintain a regulated growth-death balance. Distinguishing between these two populations is fundamental to evaluate the therapeutic ratio, as the goal of radiotherapy is to maximize tumor control while minimizing deleterious effects on healthy tissue \cite{hall2018radiobiology}.

Next, we analyze the dynamics of tumor growth in the absence of external constraints, examining three fundamental models of increasing complexity. Subsequently, we incorporate the effect of ionizing radiation through the Linear-Quadratic (LQ) model to quantify the cell lethality induced by the treatment in both populations. 

The use of the LQ model as a survival framework implies certain biological assumptions, namely that no significant repopulation or repair occurs during the irradiation itself, and that sublethal damage is fully repaired between treatment fractions. Within this context, the clinical efficacy of radiotherapy is governed by the ``5 R's'': Reoxygenation, Redistribution, Repair, Repopulation, and Radiosensitivity \cite{Steel1989}. It is worth noting that while healthy tissues repair damage more efficiently between fractions, cancerous tissues may exhibit accelerated repopulation towards the end of treatment. Furthermore, although the Gompertz model is widely used to describe these kinetics, it is limited by its prediction of a theoretical maximum volume, whereas in reality, tumor growth can often be considered indefinite.



\subsection{Free-growth}
\label{subsec:Model_No_Control_Free}

Mathematical models have long been used to describe tumor growth dynamics. The simplest formulation is the exponential model, which assumes that the growth rate is proportional to the number of tumor cells \cite{Schattler2015OptimalControlCancer}. The number of cancerous cells $C(t)$ at time $t$ is given by,
\begin{equation}
    C(t) = C_0 \, e^{\frac{\ln(2)}{T_D} (t - t_0)},
    \label{exp}
\end{equation}
where $C_0$ is the initial number of cells at time $t_0$, and $T_D$ is the constant doubling time.
This model fits well the early stages of tumor proliferation, but it overestimates growth at later stages, where nutrient depletion, hypoxia, and necrosis slow down cell division \cite{murphy2016differences}.

To address this limitation, the Gompertz model was introduced, incorporating a time-dependent decrease in proliferation rate \cite{laird1964dynamics,vaghi2020population}. Its functional form is,
\begin{equation}
     C(t) = C_0 \, e^{A \left( 1 - e^{-a (t - t_0)} \right)},
     \label{gom}
\end{equation}
where $A \in \mathbb{R}^+$ is related to the implicit carrying capacity of the tumor and $a \in \mathbb{R}^+$ controls the growth speed.

Although widely used, the Gompertz model can complicate further analysis, particularly when integrated into optimal control frameworks. It is important to note that, in mathematical oncology, the use of highly complex models with numerous parameters often does not provide additional predictive value, as the inherent uncertainty and noise in experimental and clinical data do not allow such fine considerations. Therefore, maintaining a balance between biological descriptive power and parametric simplicity is essential. In this line, we also consider the Verhulst model, which has a more abrupt growth dynamic and may fit clinical data less accurately, but offers an explicit carrying capacity \cite{forys2003logistic}. The model is given by,
\begin{equation}
    C(t) = \frac{C_0 \, e^{r t}}{1 + \frac{C_0}{K} \left( e^{r t} - 1 \right)},
    \label{verh}
\end{equation}
where $r \in \mathbb{R}^+$ is the growth rate and $K \in \mathbb{R}^+$ is the carrying capacity.

Numerical values of the parameters used in simulations are summarized in Table~\ref{tab:parameters}. These values are based on experimental data from mouse tumors, which exhibit rapid growth. Notable references include \cite{tannock1969comparison} for the doubling time $T_D$, and \cite{forys2003logistic} for the Verhulst and Gompertz parameters.

\begin{figure}[ht]
    \centering
    \includegraphics[width=0.5\linewidth]{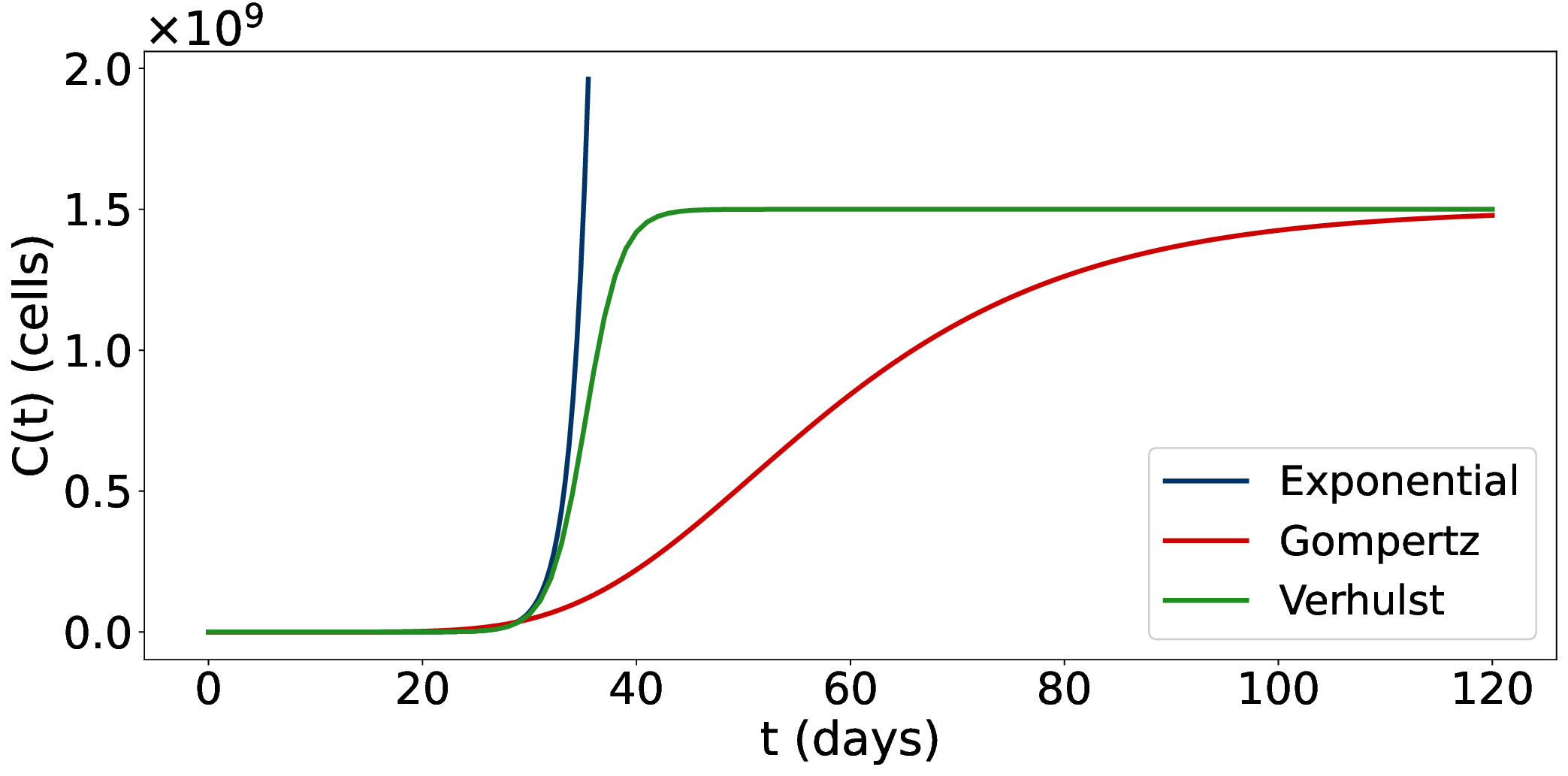}
    \caption{Numerical simulations of tumor growth dynamics in the absence of treatment, comparing the exponential, Gompertz, and Verhulst models as defined by \eqref{exp}, \eqref{gom}, and \eqref{verh}, respectively.}
    \label{fig:comparison_free}
\end{figure}

The comparative dynamics of the three growth models are illustrated in Figure~\ref{fig:comparison_free}. As expected, the exponential model \eqref{exp} fails to account for biological limitations, showing unbounded proliferation that quickly diverges from realistic clinical scales. In contrast, both the Gompertz \eqref{gom} and Verhulst \eqref{verh} models demonstrate an asymptotic convergence toward the carrying capacity $K$. Specifically, the Verhulst model exhibits a rapid, symmetric transition to saturation, while the Gompertz curve shows a more gradual deceleration in growth. This asymmetrical behavior in the Gompertz model is often more characteristic of the necrotic core formation and nutrient limitations typically observed in the development of solid tumors, \cite{alzahrani2016}.

\begin{table}[h]
	\caption{Parameter values used for the tumor growth models without treatment, corresponding to the exponential, Gompertz, and Verhulst models.}
	\centering
	{
		\begin{tabular}{ccc}
			\hline
			\makecell{\textbf{Parameter}} & 
			\makecell{\textbf{Value}} & 
			\makecell{\textbf{Description (units)}} \\
			\hline
			$C_0$ & $1$ & Initial number of cells (cells)\\
			$T_D$ & $1.15$ & Doubling time (days)  \\
			$A$ & $21.13$ & Gompertz growth rate parameter (days$^{-1}$) \\
			$a$ & $0.06$ & Gompertz shape parameter (days$^{-1}$) \\
			$K$ & $1.5 \cdot 10^9$ & Carrying capacity (cells) \\
			$r$ & $0.6$ & Verhulst growth rate (days$^{-1}$)  \\
			\hline
	\end{tabular}}
	\label{tab:parameters}
\end{table}

\subsection{Response to radiotherapy model}
\label{subsec:Model_No_Control_LQ}

To simulate the eradication of irradiated cells as a function of the dose $d$, we utilize the Linear-Quadratic (LQ) model, a cornerstone in radiobiological literature \cite{Podgorsak2005}. For a homogeneous cell population, the surviving cell number after dose $d$, the dynamics are governed by the following equation,
\begin{equation}
    N(d) = N_0 e^{-(\alpha d + \beta d^2)},
    \label{eq:LQ_model}
\end{equation}
where $d = Rt \in \mathbb{R}^+$ denotes the delivered dose and $R \in \mathbb{R}^+$ represents the dose rate. The parameters $\alpha, \beta \in \mathbb{R}^+$ characterize the two primary pathways of radiation response: lethal and sublethal damage \cite{van2018alfa}. It should be noted that the variable $N$ denotes either the cancer population $C$ or the healthy population $H$.

From a computational perspective, the system is implemented using a first-order finite difference scheme following a piecewise dynamical logic. During active radiotherapy pulses, both populations undergo a reduction phase as described by \eqref{eq:LQ_model}, where the tumor population $C$ experiences significantly greater attrition due to its inherent radiosensitivity. In the absence of radiation, the model transitions between two distinct growth regimes based on resource availability and population density. Initially, a free growth phase is assumed where both populations proliferate independently; here, the superior recovery capacity of healthy tissue is reflected by the condition $r_{free,H} > r_{free,C}$. As the tumor mass progresses and occupies the biological niche, the system enters a competition phase. In this regime, cancer cells exert dominance over the healthy tissue, following Verhulst dynamics while the healthy population remains constrained by the closed-system condition. This dominance is mathematically consistent with the assumption that $r_{free,C} > r_{comp,C}$.

\begin{figure}[htb]
    \centering
    \begin{subfigure}[b]{0.45\textwidth}
        \centering
        \includegraphics[width=\textwidth]{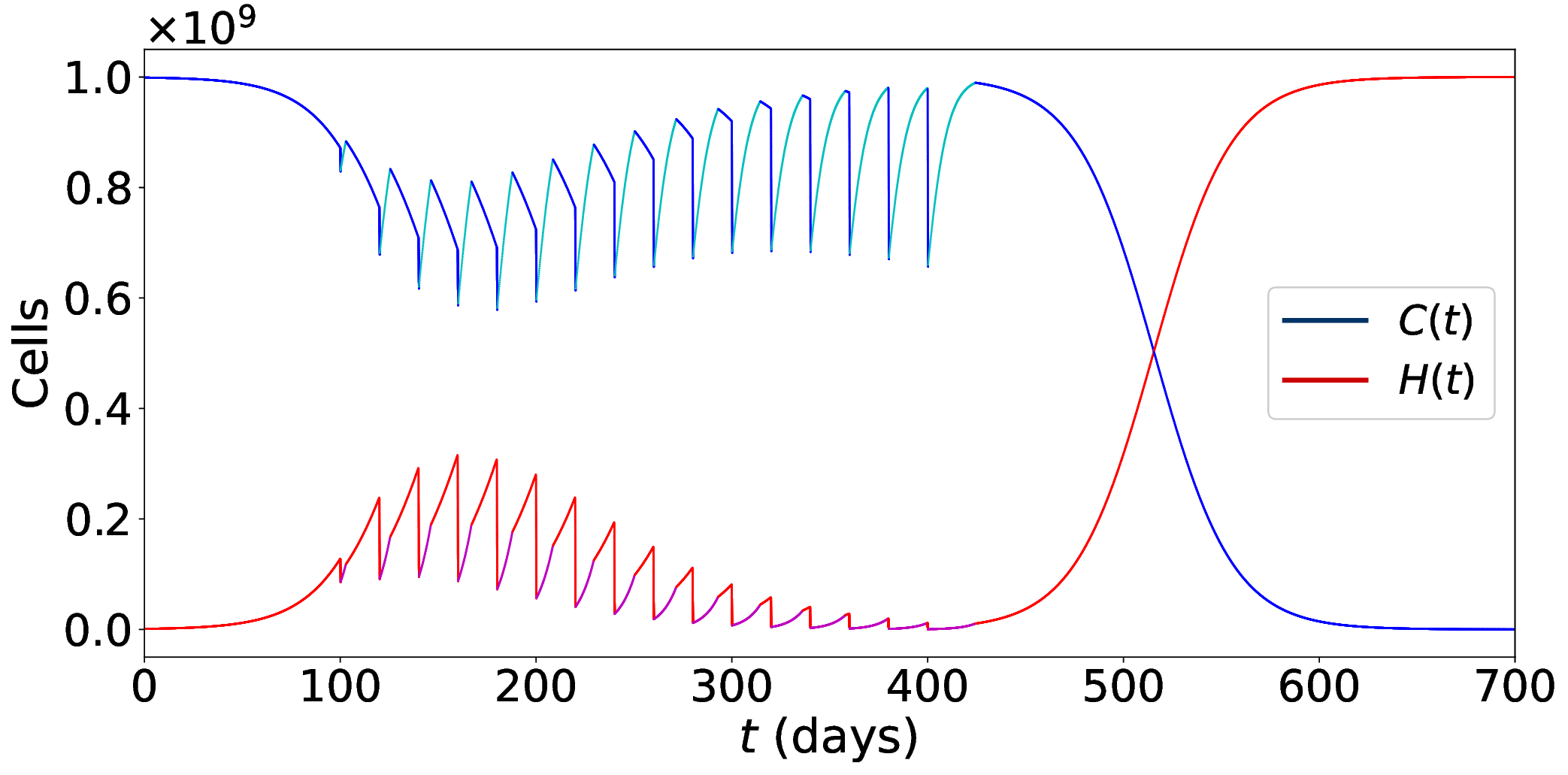}
        \caption{Cured patient (with threshold).}
        \label{fig:cured}
    \end{subfigure}
    \hfill
    \begin{subfigure}[b]{0.45\textwidth}
        \centering
        \includegraphics[width=\textwidth]{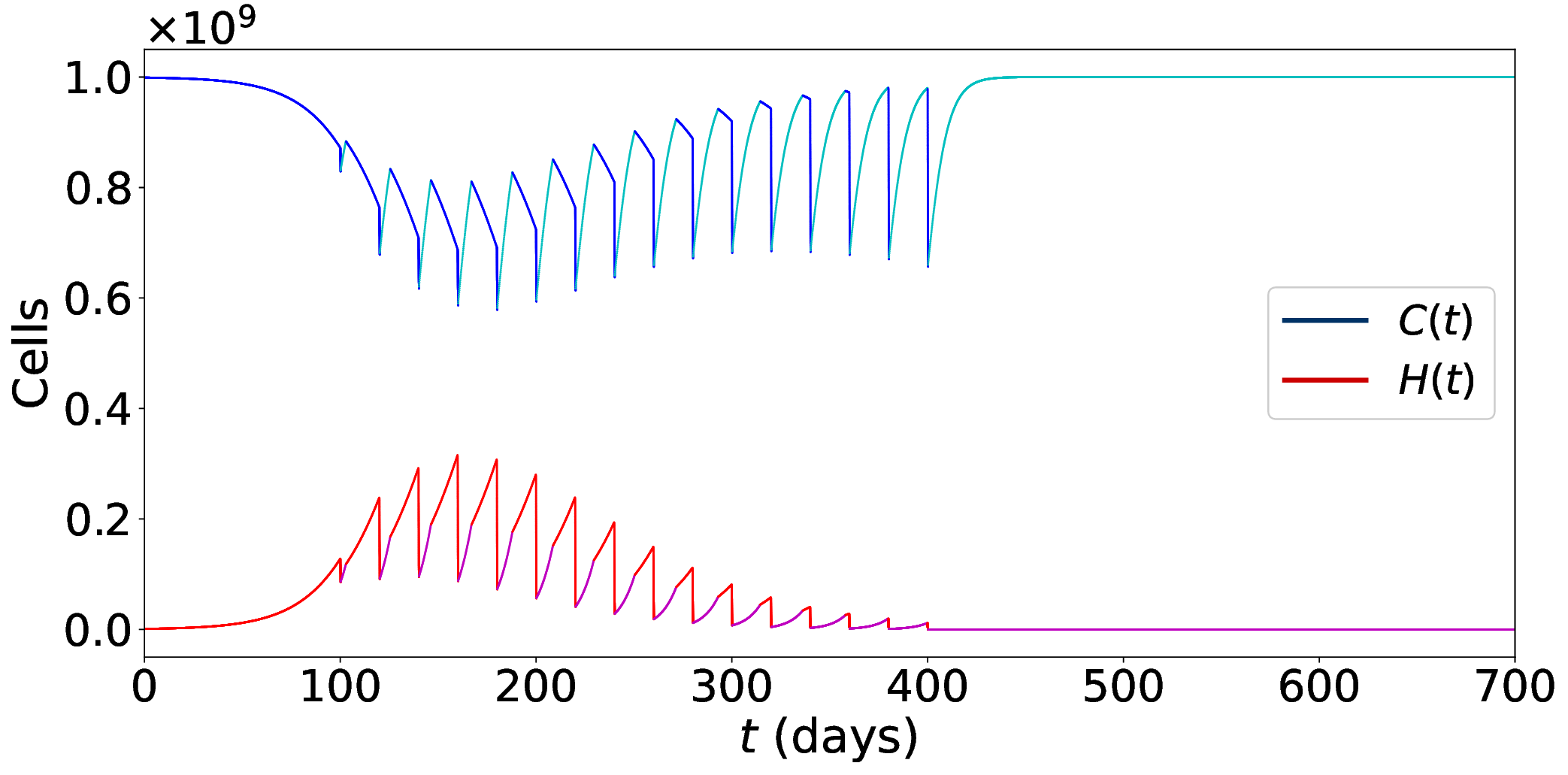}
        \caption{Non-cured patient (without threshold).}
        \label{fig:non_cured}
    \end{subfigure}
    \caption{Simulation of a fractionated radiotherapy treatment consisting of 16 sessions, each of duration $0.2$ days and separated by $20$ days. The first session is administered at $t = 100$ days and the last one at $t = 400$ days. The sharp drops correspond to treatment pulses, while the blue and red curves represent the tumor and healthy cell populations respectively, following the piecewise dynamics described in this section.}
    \label{fig:radio1D}
\end{figure}

A final consideration is the implementation of a clinical eradication threshold, $\varepsilon \in \mathbb{R}^+$. If $C(t) < \varepsilon$ following the final treatment session, the iterative process is truncated, representing a potential clinical cure. The parameters utilized for the simulation in Figure~\ref{fig:radio1D} are summarized in Table \ref{tab:radiotherapy_parameters}. To ensure clinical relevance, we adopted a human-like growth rate, $r < 0.6$ days$^{-1}$, as proposed in \cite{tannock1969comparison}. Regarding the radiotherapy parameters, we acknowledge that the dose rate $R = 1$ Gy/day and the session duration of 0.2 days are simplified numerical approximations. In clinical practice, dose rates are significantly higher ($R \sim 1$ Gy/min), delivered in shorter intervals. However, to facilitate the visualization of the dynamical transitions and the competitive interaction between populations over a 700-day period, we have scaled these parameters accordingly. This ensures the qualitative behavior of the system remains numerically observable within the ODE framework while maintaining biologically plausible $\alpha/\beta$ ratios between $1$ and $10$ Gy \cite{hall2018radiobiology}. To compensate for this temporal scaling, the parameters $\alpha$ and $\beta$ were adjusted to reflect a more aggressive therapeutic regimen that remains consistent with the differential radiosensitivity of the tissues.

As shown in Figure~\ref{fig:cured}, the tumor population (blue curve) falls below the clinical eradication threshold after the final treatment session, allowing the healthy tissue (red curve) to recover the niche. In contrast, in Figure~\ref{fig:non_cured}, the tumor regrows during the inter-fraction intervals, ultimately leading to treatment failure.

\begin{table}[ht]
    \caption{Parameter values used for the radiotherapy simulation based on the piecewise Verhulst-LQ model. Note: $R$, $\alpha$, and $\beta$ are scaled for numerical visualization while preserving clinical $\alpha/\beta$ ratios.}
    \centering
    {
        \begin{tabular}{ccc}
            \hline
            \makecell{\textbf{Parameter}} & 
            \makecell{\textbf{Value}} & 
            \makecell{\textbf{Description (units)}} \\
            \hline
            $C_0$ & $10^6$ & Initial cancer cell population (cells) \\
            $K$ & $10^9$ & Carrying capacity (cells) \\
            $\varepsilon$ & $10^6$ & Non-clonogenic threshold (cells) \\
            $R$ & $1$ & Radiation dose rate (Gy/day) \\
            $r_{free,C}$ & $0.13$ & Growth rate, cancer cells (days$^{-1}$) \\
            $r_{comp,C}$ & $0.05$ & Competitive growth rate, cancer cells (days$^{-1}$) \\
            $r_{free,H}$ & $0.16$ & Growth rate, healthy cells (days$^{-1}$) \\
            $\alpha_C$ & $5 \cdot 10^{-3}$ & Linear inactivation coefficient, cancer (Gy$^{-1}$) \\
            $\beta_C$ & $2 \cdot 10^{-2}$ & Quadratic inactivation coefficient, cancer (Gy$^{-2}$) \\
            $\alpha_H$ & $6.25 \cdot 10^{-4}$ & Linear inactivation coefficient, healthy (Gy$^{-1}$) \\
            $\beta_H$ & $2.5 \cdot 10^{-3}$ & Quadratic inactivation coefficient, healthy (Gy$^{-2}$) \\
            \hline
    \end{tabular}}
    \label{tab:radiotherapy_parameters}
\end{table}

\subsection{Advanced model of tumor dynamics}
\label{subsec:Model_No_Control_Advanced}

In previous sections, tumor growth was modeled under the assumption of a closed-system condition. This simplification allowed for the characterization of the system by tracking a single cell population, deriving the other as a dependent variable. However, clinical evidence increasingly suggests that tumor progression is inextricably linked to complex chemical and energetic exchanges with the host microenvironment. For instance, \cite{lee2015crosstalk} demonstrated the bidirectional information exchange between tumor cells and the hematological or lymphatic systems, particularly during tumor-induced angiogenesis. Similar findings regarding revascularization processes are discussed by \cite{emami2021role}.

To better reflect these biological realities, we transition toward a more flexible framework by relaxing the closed-system constraint. Following a mathematical approach similar to the epidemiological study by \cite{lopez2024optimal}, we propose a system of two coupled ordinary differential equations,
\begin{equation}
\begin{split}
    \dfrac{dH(t)}{dt} &= r_H \left( 1 - \dfrac{H(t) + C(t)}{K_H} \right) H(t) \\[8pt]
    \dfrac{dC(t)}{dt} &= r_C \left( 1 - \dfrac{H(t) + C(t)}{K_C} \right) C(t),
\end{split}
\label{eq:coexistence_system}
\end{equation}
where $H(t)$ and $C(t)$ represent the healthy and cancerous cell populations, respectively. The parameters $r_H, r_C \in \mathbb{R}^+$ denote their intrinsic growth rates, while $K_H, K_C \in \mathbb{R}^+$ represent their respective carrying capacities.

Two fundamental aspects of this model must be emphasized. First, since our simulation focuses on the treatment window, we assume $r_H > r_C$.  As is well established in clinical oncology, radiotherapy is typically delivered in a fractionated and moderated manner in order to minimize damage to healthy tissues while maximizing its effect on tumor cells. This strategy allows healthy tissues to regenerate between treatment sessions, whereas the tumor is progressively eliminated, since clonogenic cancer cells progressively lose their reproductive capacity. In this context, healthy tissues exhibit superior regenerative capacity. In our previous uncoupled models, the dominant nature of the tumor could only be enforced by assuming $r_H < r_C$ under the closed-system condition. However, the current coupled formulation allows us to maintain $r_H > r_C$ while accounting for the competitive advantage of malignant cells observed after treatment sessions.

Second, we must reconsider the ecological interpretation of the carrying capacity $K$, defined as the maximum population size that an environment can sustain without degradation. While it might be argued that $K_C > K_H$, given that cancer cells often bypass homeostatic constraints, such anomalous growth ultimately leads to the systemic deterioration of the host, which contradicts the sustainability requirement of the carrying capacity. Consequently, one might hypothesize that $K_C = K_H = K$. However, this assumption leads to a mathematical contradiction: even with distinct growth rates $r$, a system where $K_C = K_H$ tends toward coexistence rather than the competitive exclusion or dominance characteristic of malignant progression.

Mathematically, the system's behavior can be elucidated by solving the associated eigenvalue problem. Following the methodology described in \cite{strogatz2018nonlinear} and focusing on the primary points of interest, we obtain the results presented in Table \ref{tab:stability_analysis}.

\begin{table}[ht]
\centering
\renewcommand{\arraystretch}{1.5} 
\caption{Equilibrium points and stability analysis for the coexistence model described by system \eqref{eq:coexistence_system}.}
\label{tab:stability_analysis}
\begin{tabular}{ccc}
\hline
\textbf{Equilibrium Point} & \textbf{Eigenvalues} & \textbf{Stability} \\ \hline
\multirow{2}{*}{$(0, 0)$} & $\lambda_1 = r_H > 0$ & \multirow{2}{*}{Unstable} \\
                           & $\lambda_2 = r_C > 0$ & \\ \hline
\multirow{2}{*}{$(K_H, 0)$} & $\lambda_1 = -r_H < 0$ & \multirow{2}{*}{Unstable} \\
                             & $\lambda_2 = \frac{r_C(K_C-K_H)}{K_C} = 0$ & \\ \hline
\multirow{2}{*}{$(0, K_C)$} & $\lambda_1 = -r_C < 0$ & \multirow{2}{*}{Unstable} \\
                             & $\lambda_2 = \frac{r_H(-K_C+K_H)}{K_H} = 0$ & \\ \hline
\end{tabular}
\end{table}

\begin{figure}[hbt]
    \centering
    \includegraphics[width=0.75\textwidth]{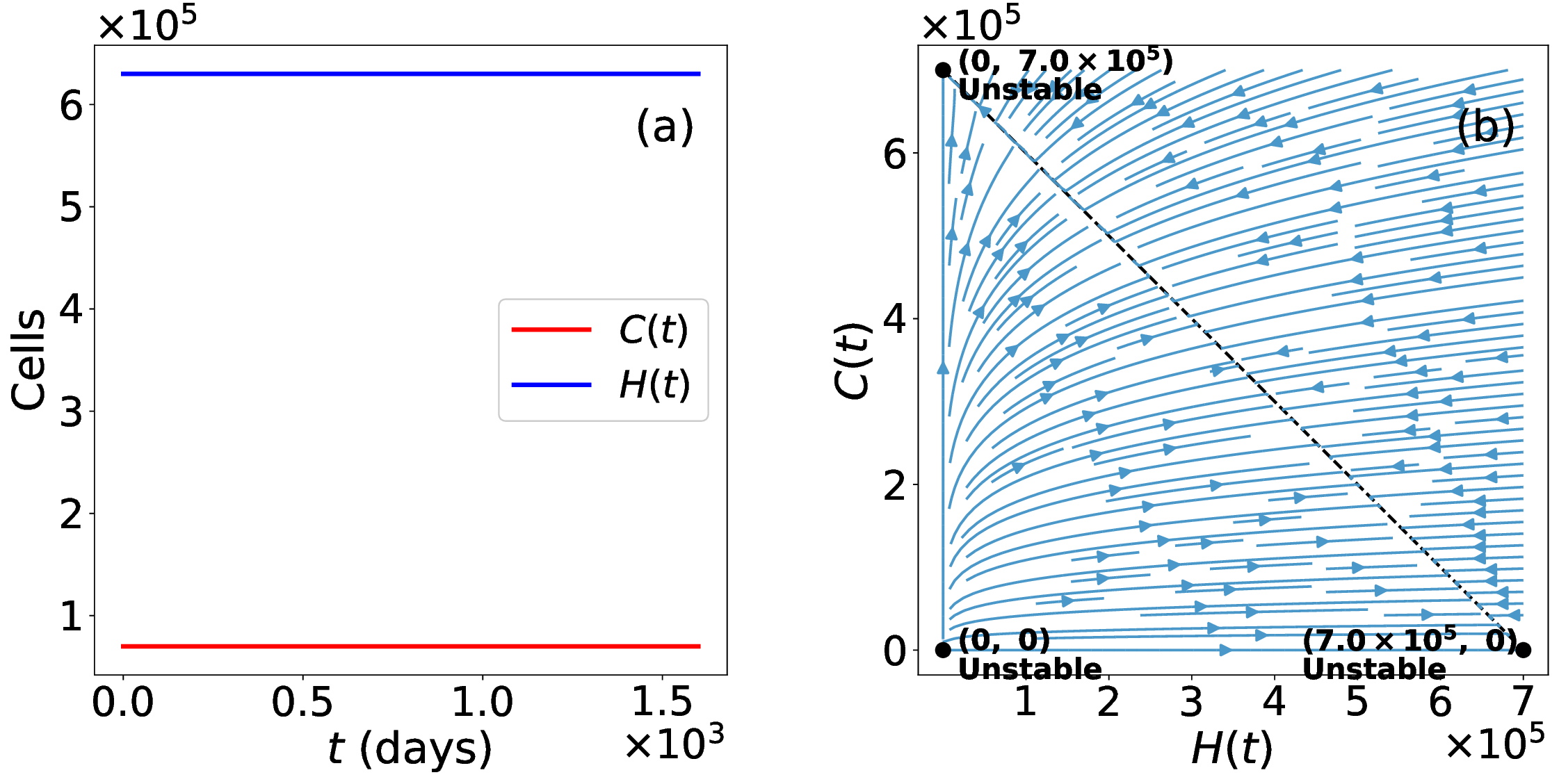}
    \caption{(a) Temporal evolution of healthy and cancer cell populations in the coexistence mechanisms. (b) Phase portrait of system \eqref{eq:coexistence_system}, illustrating the vector field and the instability of the equilibrium points.}
    \label{fig:placeholder}
\end{figure}

Graphically, these dynamical features are illustrated through the temporal evolution and phase portraits shown in Figure~\ref{fig:placeholder}. For these simulations, we utilize the numerical values summarized in Table~\ref{tab:model_parameters}, which are based on empirical data from murine tumor dynamics \cite{forys2003logistic}. As specified in the parameters, we assume the healthy tissue regenerates five times faster than the cancerous cells, $r_H = 5r_C$.

The results in Figure~\ref{fig:placeholder}a confirm that, despite this significant regenerative advantage of the healthy cells, the system does not settle into a stable state. Instead, we observe a transient regime where both populations coexist without reaching a definitive ecological climax. This behavior is further elucidated by the vector field in the phase portrait (Figure~\ref{fig:placeholder}b), which reflects the lack of stable attractors identified in Table~\ref{tab:stability_analysis}. The trajectories are repelled from the axis, particularly from the healthy equilibrium $(K, 0)$, confirming that the system is topologically unable to return to a tumor-free state autonomously.

Nonetheless, to ensure numerical stability and convergence within the optimization framework described in Section~\ref{sec:Model_Control}, a scaled value for $K$ was adopted. Large differences in the order of magnitude between state variables and their sensitivities often lead to ill-conditioned gradients; therefore, reducing the carrying capacity facilitates a more robust numerical resolution while preserving the qualitative dynamics of the system. This allows for a coherent comparison of the results obtained across different contexts.

\begin{table}[ht]
\centering
\caption{Parameter values and initial conditions for the coexistence model.}
\label{tab:model_parameters}
\begin{tabular}{cc}
\hline
\textbf{Parameter} & \textbf{Value (units)} \\ \hline
$r_C$ & $0.6$ (days$^{-1}$) \\
$r_H$ & $3$ (days$^{-1}$) \\
$K$ & $7 \cdot 10^5$ (cells) \\
$C_0$ & $0.7 \cdot 10^5$ (cells) \\
$H_0$ & $6.3 \cdot 10^5$ (cells) \\ \hline
\end{tabular}
\end{table}

Indeed, as anticipated, no inter-species competition is observed in the previous formulation; instead, a coexistence regime emerges, reflected by the existence of a line of unstable equilibria in the phase diagram.

Consequently, for our model to be valid for describing a tumor system, we must introduce a small competition term that favors the cancer cells in system \eqref{eq:coexistence_system}, thus orienting the stability toward the point $(0, K)$. In this way, our system of equations is redefined as follows,
\begin{equation}
\begin{split}
    \dfrac{dH(t)}{dt} &= r_H \left( 1 - \dfrac{H(t) + C(t)}{K} \right) H(t) - \gamma H(t)C(t) \\[8pt]
    \dfrac{dC(t)}{dt} &= r_C \left( 1 - \dfrac{H(t) + C(t)}{K} \right) C(t),
\end{split}
\label{eq:competition_system}
\end{equation}
where $\gamma = 5.5 \times 10^{-8}$ days$^{-1}$ cells$^{-1}$ is the factor reflecting the magnitude of the competition. The new equilibrium points and their stability are summarized in Table \ref{tab:stability_competition}.

\begin{table}[ht]
\centering
\renewcommand{\arraystretch}{1.8}
\caption{Equilibrium points and stability for the competition model without treatment, given by system \eqref{eq:competition_system}.}
\label{tab:stability_competition}
\begin{tabular}{ccc}
\hline
\textbf{Equilibrium Point} & \textbf{Eigenvalues} & \textbf{Stability} \\ \hline
\multirow{2}{*}{$(0, 0)$} & $\lambda_1 = r_H > 0$ & \multirow{2}{*}{Unstable} \\
                           & $\lambda_2 = r_C > 0$ & \\ \hline
\multirow{2}{*}{$(K, 0)$} & $\lambda_1 = -r_H < 0$ & \multirow{2}{*}{Unstable} \\
                           & $\lambda_2 = r_C \frac{K-K}{K} = 0$ & \\ \hline
\multirow{2}{*}{$(0, K)$} & $\lambda_1 = -r_C < 0$ & \multirow{2}{*}{Stable} \\
                           & $\lambda_2 = r_H \frac{-K+K}{K} - \gamma K = -\gamma K < 0$ & \\ \hline
\end{tabular}
\end{table}

\begin{figure}
    \centering
    \includegraphics[width=0.75\textwidth]{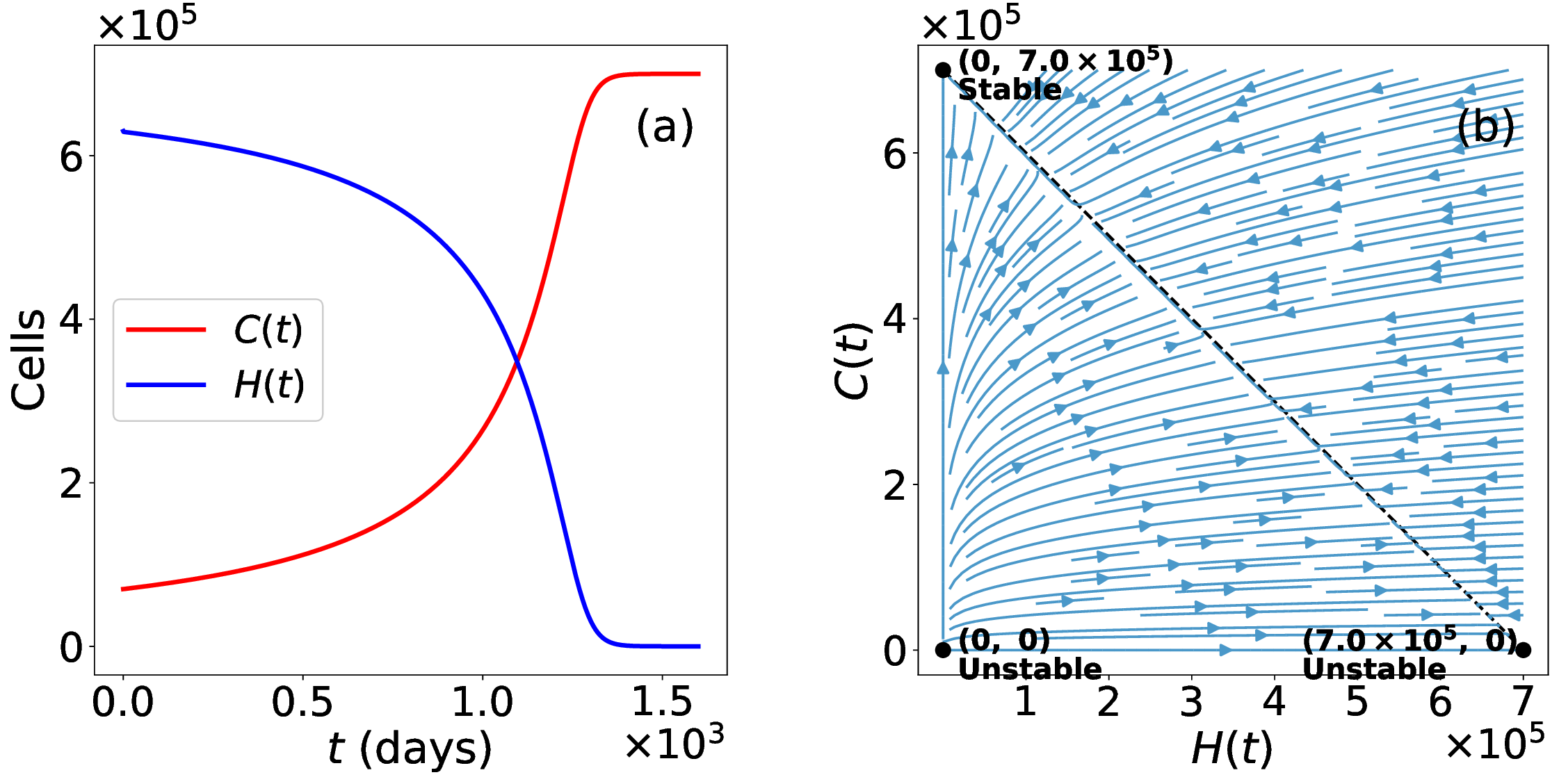}
    \caption{(a) Temporal evolution of healthy $H(t)$ and cancer $C(t)$ cell populations under the competition model \eqref{eq:competition_system}. (b) Phase portrait of the autonomous system showing the vector field. The labels indicate the transition of $(0, K)$ into a stable sink and $(K, 0)$ into an unstable saddle point, as derived in Table \ref{tab:stability_competition}.}
    \label{fig:placeholder2}
\end{figure}

Graphically, the new situation is reflected in the phase portrait and temporal evolution, where the previous coexistence regime is replaced by the clear dominance of the malignant population. As observed in Figure~\ref{fig:placeholder2}a, even with the regeneration rate of healthy cells being significantly higher, $r_H = 5r_C$, the inclusion of the competition term $\gamma$ causes a steady decline in $H(t)$ as the tumor progresses toward the carrying capacity $K$. This shift is mathematically confirmed in the phase portrai, Figure~\ref{fig:placeholder2}b, where the point $(0, K)$ now acts as a stable sink, while the healthy equilibrium $(K, 0)$ has become an unstable saddle.

\begin{figure}
    \centering
    \includegraphics[width=\textwidth]{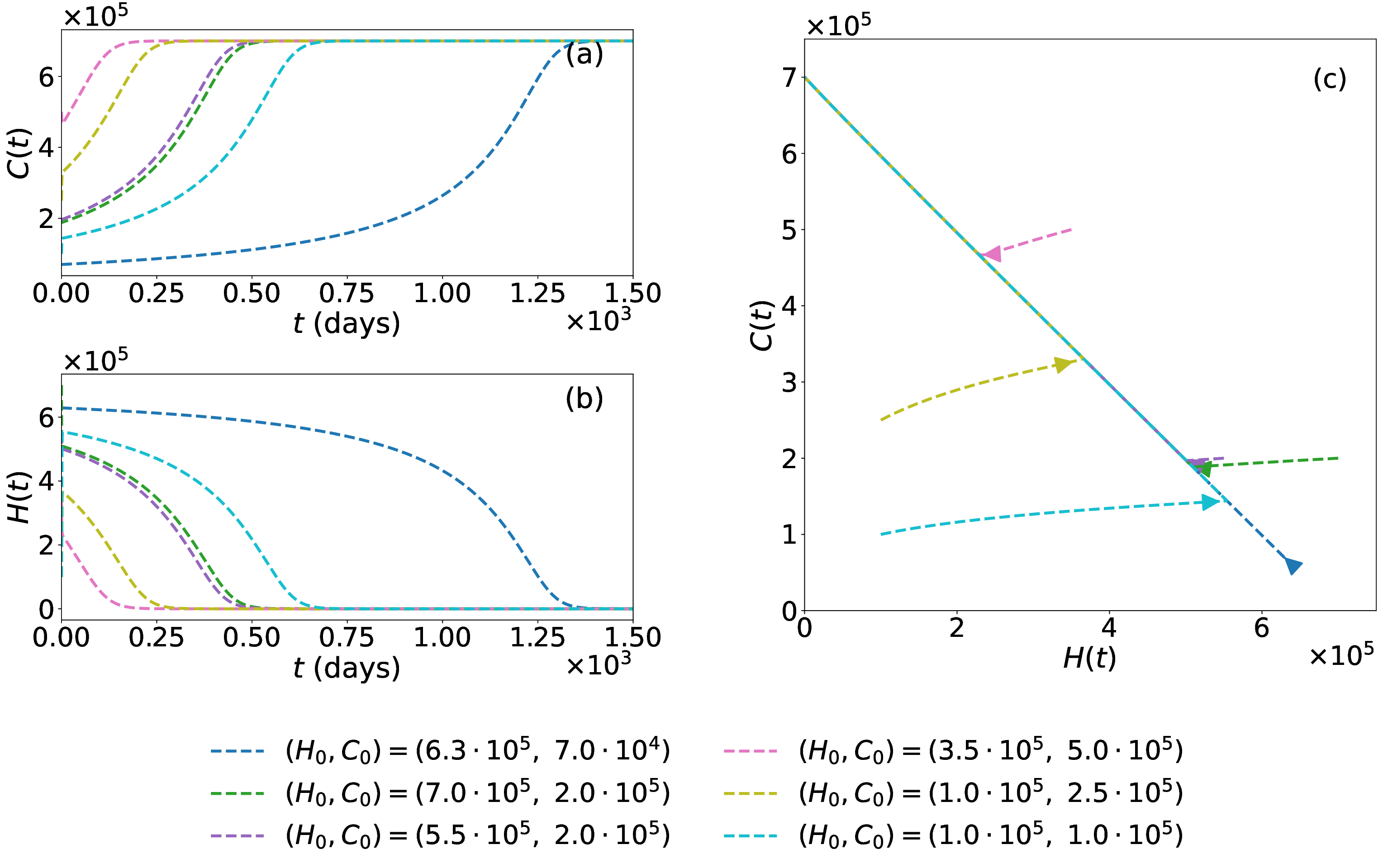}
    \caption{Sensitivity of the competition model to initial conditions using parameters from Table \ref{tab:model_parameters}. (a) Temporal progression of the tumor population $C(t)$ and (b) decline of the healthy population $H(t)$ for various $(H_0, C_0)$ pairs. (c) Combined phase portrait illustrating that all trajectories, regardless of the starting density, are globally attracted to the tumor-dominant equilibrium $(0, K)$.}
    \label{fig:placeholder3}
\end{figure}

The robustness of this competitive advantage is further demonstrated in Figure~\ref{fig:placeholder3}. Regardless of the initial conditions $(H_0, C_0)$ chosen from the murine-based values in Table~\ref{tab:model_parameters}, the system's trajectories consistently converge toward the tumor-only state. Figures~\ref{fig:placeholder3}a and \ref{fig:placeholder3}b show that even when starting with a high density of healthy cells and a minimal tumor presence, the competitive pressure $\gamma = 5.5 \times 10^{-8}$ eventually suppresses the healthy tissue. 

The phase portrait in Figure~\ref{fig:placeholder3}c summarizes this behavior: the vector field no longer supports a line of equilibria but instead directs all trajectories toward $(0, K)$. This topological change confirms that without external intervention, the tumor is guaranteed to displace the healthy population, providing a definitive biological and mathematical justification for the optimal control problem addressed in the following section.

	\section{Tumor model with control}
	\label{sec:Model_Control}
	
In the field of mathematical oncology, optimal control theory has become a cornerstone for designing treatment protocols that balance tumor regression with the preservation of healthy tissue. Unlike empirical protocols, this mathematical framework allows for the minimization of the total tumor burden while strictly limiting the toxicity induced by the intervention \cite{Schattler2015OptimalControlCancer}. In this section, we apply these principles to invert the stability landscape of the cellular competition model previously shown. Our objective is to transition the system from a state of malignant dominance toward a stable healthy equilibrium $(K, 0)$. This is achieved by modulating the radiotherapy intensity through a control parameter $u(t)$, which represents the dose rate. By employing optimal control techniques, we aim to minimize $u(t)$ to ensure that the therapeutic gain is maximized while collateral damage to the biological niche is kept at a minimum.

	\subsection{Radiotherapy optimal control problem}

From a mathematical perspective, the inclusion of radiotherapy involves the addition of new terms to system \eqref{eq:competition_system} that reflect damage proportional to the population density, while maintaining a higher aggression toward cancerous cells. 

Following the logic applied in multi-strain viral models, where treatment efficacy varies between lineages, we introduce the constants $\lambda$ and $\mu$ to regulate the impact of the therapy. While in virology a treatment might be more effective against an original strain ($c_A > c_B$), in the context of oncology, we define $\lambda < \mu$ to account for the greater attrition experienced by the tumor population due to its radiosensitivity. The dynamics of the controlled system are thus given by,
\begin{equation}
\begin{split}
    \dfrac{dH(t)}{dt} &= r_H \left( 1 - \dfrac{H(t) + C(t)}{K} \right) H(t) - \gamma H(t)C(t) - \lambda u(t) H(t) \\[8pt]
    \dfrac{dC(t)}{dt} &= r_C \left( 1 - \dfrac{H(t) + C(t)}{K} \right) C(t) - \mu u(t) C(t),
\end{split}
\label{eq:controlled_system}
\end{equation}
where $r_H, r_C, K, \gamma \in \mathbb{R}^+$, $\lambda, \mu \in [0, 1]$, additionally $\lambda < \mu$ and the control function $u(\cdot)$ satisfies the constraint $0 \leq u(\cdot) \leq u_{max} \leq 1$.

The set of admissible control functions is given by,
\begin{equation}
    \Omega = \{u(\cdot) \in L^\infty(0, t_f) \mid 0 \leq u(t) \leq u_{max} \leq 1 , \forall t \in [0, t_f] \} .
\end{equation}

When $u(t) = 0$, no treatment is supplied to the system; i.e., we recover the original competition model proposed in \eqref{eq:competition_system}. On the other hand, when $u(t) \in (0, 1]$, we are introducing radiotherapy to the system, which reduces the cellular populations. Moreover, we consider constant initial conditions $(H(0), C(0)) = (H_0, C_0) \in \mathbb{R}^+ \times \mathbb{R}^+$.

The main goal is to minimize the tumor population $C(t)$ in the time window $[0, t_f]$, effectively forcing the system toward the healthy fixed point $(K, 0)$ with the least possible cost. Therefore, following \cite{ibanez2017optimal}, we define the following cost functional, $J$:

\begin{equation}
    J(H(t), C(t), u(t)) = \int_{0}^{t_f} \left[ (H(t) - K)^2 + (C(t) - 0)^2 + u^2(t) \right] dt.
    \label{eq:cost_functional}
\end{equation}

This functional is defined to drive the healthy cell density $H$ toward its carrying capacity $K$ and the cancer population $C$ toward extinction. In other words, we seek to minimize the cost of approaching the healthy-dominant state while minimizing the total radiation dose to prevent excessive collateral damage. The existence of solutions for the optimal control problem \eqref{eq:controlled_system}-\eqref{eq:cost_functional} is ensured by classical sufficient conditions, see e.g. \cite{trelat2005controle} and references cited therein.

\subsubsection{\texorpdfstring{Analytical study of the system~\eqref{eq:controlled_system} with constant radiotherapy treatment}{Analytical study of the controlled system with constant radiotherapy treatment}}

To evaluate the impact of a continuous therapeutic strategy, we first consider the case where the control parameter is constant, $u(t) = {u} \in (0, 1]$. This scenario represents a simplified clinical protocol where radiotherapy is administered at a uniform intensity throughout the observation period. Under this assumption, the system of differential equations \eqref{eq:controlled_system} is simplified, and the dynamics are governed by,
\begin{equation}
\begin{split}
    \dfrac{dH(t)}{dt} &= r_H \left( 1 - \dfrac{H(t) + C(t)}{K} \right) H(t) - \gamma H(t)C(t) - \lambda {u} H(t) \\[8pt]
    \dfrac{dC(t)}{dt} &= r_C \left( 1 - \dfrac{H(t) + C(t)}{K} \right) C(t) - \mu {u} C(t),
\end{split}
\label{eq:constant_control}
\end{equation}
where ${u}$ is the constant dosage level. The inclusion of this term allows for a direct analysis of how the stability of the equilibrium points shifts as a function of the treatment intensity.

The equilibrium points in the  $H$-$C$ plane result modified once the radiotherapy treatment is introduced, represented by the fixed treatment rate $u$. In this subsection we study the behavior of the system if we apply a constant dose of radiotherapy over time, i.e., when a constant value is considered for the control, $u$. 
The fixed points are obtained by setting the derivatives to zero, which is equivalent to solving the following algebraic system for the non-trivial cases,
\begin{equation}
\begin{split}
    r_H \left( 1 - \dfrac{H(t) + C(t)}{K} \right) - \gamma C(t) - \lambda {u} = 0 \\
    r_C \left( 1 - \dfrac{H(t) + C(t)}{K} \right) - \mu {u} = 0.
\end{split}
\end{equation}

The stability of these equilibrium points is determined by the Jacobian matrix of the system,
\begin{equation}
J(H(t), C(t)) = \begin{pmatrix}
r_H \left( 1 - \frac{2H(t) + C(t)}{K} \right) - \gamma C(t) - \lambda {u} & -H(t) \left( \frac{r_H}{K} + \gamma \right) \\
-C(t) \frac{r_C}{K} & r_C \left( 1 - \frac{H(t) + 2C(t)}{K} \right) - \mu {u}
\end{pmatrix}.
\end{equation}

By evaluating the eigenvalues of this matrix at each fixed point, we can characterize the stability landscape as a function of the treatment intensity ${u}$. The results of this stability analysis are summarized in Table~\ref{cond_control}.

\begin{table}[ht]
\centering
\renewcommand{\arraystretch}{1.8}
\caption{Equilibrium points and stability for the competition model with constant treatment $u$, given by system \eqref{eq:constant_control}.}
\label{cond_control}
\resizebox{\textwidth}{!}{
\begin{tabular}{ccc}
\hline
\textbf{Equilibrium Point} & \textbf{Eigenvalues} & \textbf{Stability} \\ \hline
\multirow{2}{*}{$(0, 0)$} & $\lambda_1 = r_H - \lambda u$ & \multirow{2}{*}{$u > \max\left\{\frac{r_H}{\lambda}, \frac{r_C}{\mu} \right\}$} \\
                           & $\lambda_2 = r_C - \mu u$ & \\ \hline
\multirow{2}{*}{$\left(K \left( 1 - \frac{\lambda u}{r_H} \right), 0\right)$} & $\lambda_1 = \lambda u - r_H$ & \multirow{2}{*}{\makecell{$u < \frac{r_H}{\lambda}$ \\ $\frac{r_C}{r_H} < \frac{\mu}{\lambda}$}} \\
                           & $\lambda_2 = u \left( \frac{r_C \lambda - r_H \mu}{r_H} \right)$ & \\ \hline
\multirow{2}{*}{$\left(0, K \left( 1 - \frac{\mu u}{r_C} \right)\right)$} &  $\lambda_1 = \mu u - r_C$  & \multirow{2}{*}{\makecell{$u < \frac{r_C}{\mu}$ \\ $\gamma > \frac{u}{K} \left( \frac{r_H \mu - \lambda r_C}{r_C - \mu u} \right)$ }} \\
                           & $\lambda_2 = \frac{u(r_H \mu - \lambda r_C)}{r_C} - \gamma K \left( 1 - \frac{\mu u}{r_C} \right)$& \\ \hline
\multirow{2}{*}{\makecell{$(\hat{H}, \hat{C}) = $ \\ $\left( K - \frac{K \mu u}{r_C} + \frac{\lambda u}{\gamma} - \frac{\mu r_H u}{\gamma r_C} , \frac{u(\mu r_H - \lambda r_C)}{\gamma r_C} \right)$}} & \multirow{2}{*}{-} & \multirow{2}{*}{\makecell{Unstable \\ $\left(\det(J(\hat{H},\hat{C})<0)\right)$}} \\
                           & & \\ \hline
\end{tabular}}
\end{table}

The analytical framework developed above provides a rigorous characterization of the system's equilibria under a constant radiotherapy regimen. Specifically, the conditions derived in Table \ref{cond_control} define the threshold values for the treatment intensity $u$ necessary to alter the competitive balance between healthy and malignant cells, potentially leading to tumor eradication.

\begin{table}[ht]
\centering
\renewcommand{\arraystretch}{1.5}
\caption{Summary of parameter values for the competition model with control, as defined by system \eqref{eq:controlled_system}.}
\label{tab:control_parameters}
\begin{tabular}{cc}
\hline
\textbf{Parameter} & \textbf{Value (units)} \\ \hline
$r_C$ & $0.6$ (days$^{-1}$) \\
$r_H$ & $3$ (days$^{-1}$) \\
$K$ & $7 \times 10^5$ (cells) \\
$\gamma$ & $5.5 \times 10^{-8}$ (cells$^{-1}$ days$^{-1}$) \\
$\lambda$ & $0.025$ (Gy$^{-1}$) \\
$\mu$ & $0.189$ (Gy$^{-1}$) \\ \hline
\end{tabular}
\end{table}

\begin{figure}[ht]
		\centering
		\includegraphics[width=\textwidth]{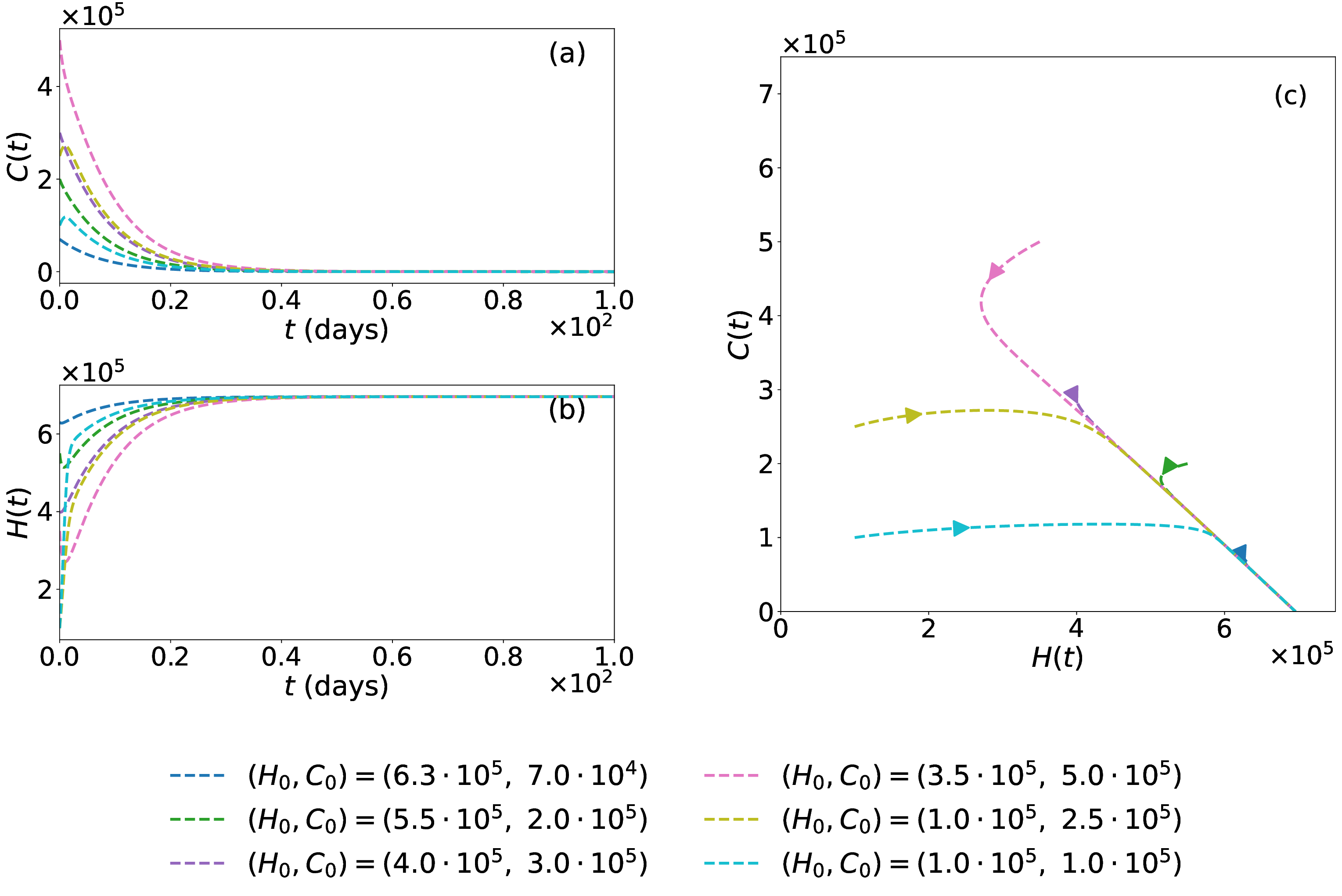}
        \caption{System dynamics under a constant radiotherapy treatment $u=0.7\, Gy$. Results are shown for various initial conditions $(H_0, C_0)$ using the parameter set from Table \ref{tab:control_parameters}. (a) Temporal evolution of $C(t)$ showing eradication. (b) Temporal evolution of $H(t)$ illustrating recovery toward the shifted equilibrium. (c) Phase portrait $H(t)$ vs. $C(t)$, where the vector field demonstrates that the therapeutic pressure $u$ makes the point $\left(K \left( 1 - \frac{\lambda u}{r_H} \right), 0\right)$ a globally stable attractor.}
		\label{fig:model_constant_control}
\end{figure}

The numerical results presented in Figure \ref{fig:model_constant_control} illustrate the effectiveness of the constant treatment $u$ in fulfilling the stability conditions derived in Table \ref{cond_control}. For these simulations, a constant dose is applied, satisfying the threshold required to destabilize the tumor-dominant state and promote the recovery of the healthy population.

As shown in Figures \ref{fig:model_constant_control}a and \ref{fig:model_constant_control}b, regardless of the initial tumor burden, the cancerous population $C(t)$ exhibits an exponential decay toward zero, while the healthy population $H(t)$ recovers. However, as predicted by the analytical study, the healthy tissue does not return to its original carrying capacity $K$, but settles at the new equilibrium level $K \left( 1 - \frac{\lambda u}{r_H} \right)$ due to the toxicity of the continuous treatment.

The phase portrait in Figure \ref{fig:model_constant_control}c provides a global perspective of this transition. All trajectories, representing diverse initial conditions, are attracted to this new healthy-fixed point. This confirms that the constant treatment $u$ has successfully transformed the previously unstable healthy region into a stable sink, effectively collapsing the tumor's biological niche at the cost of a controlled reduction in the total healthy cell count.

\subsubsection{Numerical simulations for the optimal control problem using a direct method}

Numerical simulations were performed using the \textit{OptimalControl.jl} toolbox \cite{Caillau_OptimalControl_jl_a_Julia} implemented in the \textit{Julia} programming language.

In order to implement the problem in \textit{Julia}, the optimal control problem is reformulated by introducing an auxiliary state variable, $z(t)$, associated with the cost function. This transformation allows the problem to be expressed in \textit{Mayer form}, where the cost functional is defined solely in terms of the final state of the augmented system \cite{Liberzon2012}. The resulting augmented system is,
\begin{equation}
    \tag{OCP}
    \begin{cases}
        z(t_f) + \displaystyle\int_0^{t_f} u^2(t) \ dt \rightarrow \min \\[0.2cm]
        \dot{H}(t) = r_H \left( 1 - \dfrac{H(t) + C(t)}{K} \right) H(t) - \gamma H(t)C(t) - \lambda u(t) H(t) \\
        \dot{C}(t) = r_C \left( 1 - \dfrac{H(t) + C(t)}{K} \right) C(t) - \mu u(t) C(t) \\
        \dot{z}(t) = (H(t) - K)^2 + C(t)^2
        \\[0.2cm]
        u(t) \in [0, u_{\max}], \ t \in [0, t_f] \\
        H(t), C(t), z(t) \geq 0 \\
        H(0) = H_0, \ C(0) = C_0, \ z(0) = 0.
    \end{cases}
    \label{OCP}
\end{equation}
	
	\begin{figure}[ht]
		\centering
		\includegraphics[width=\textwidth]{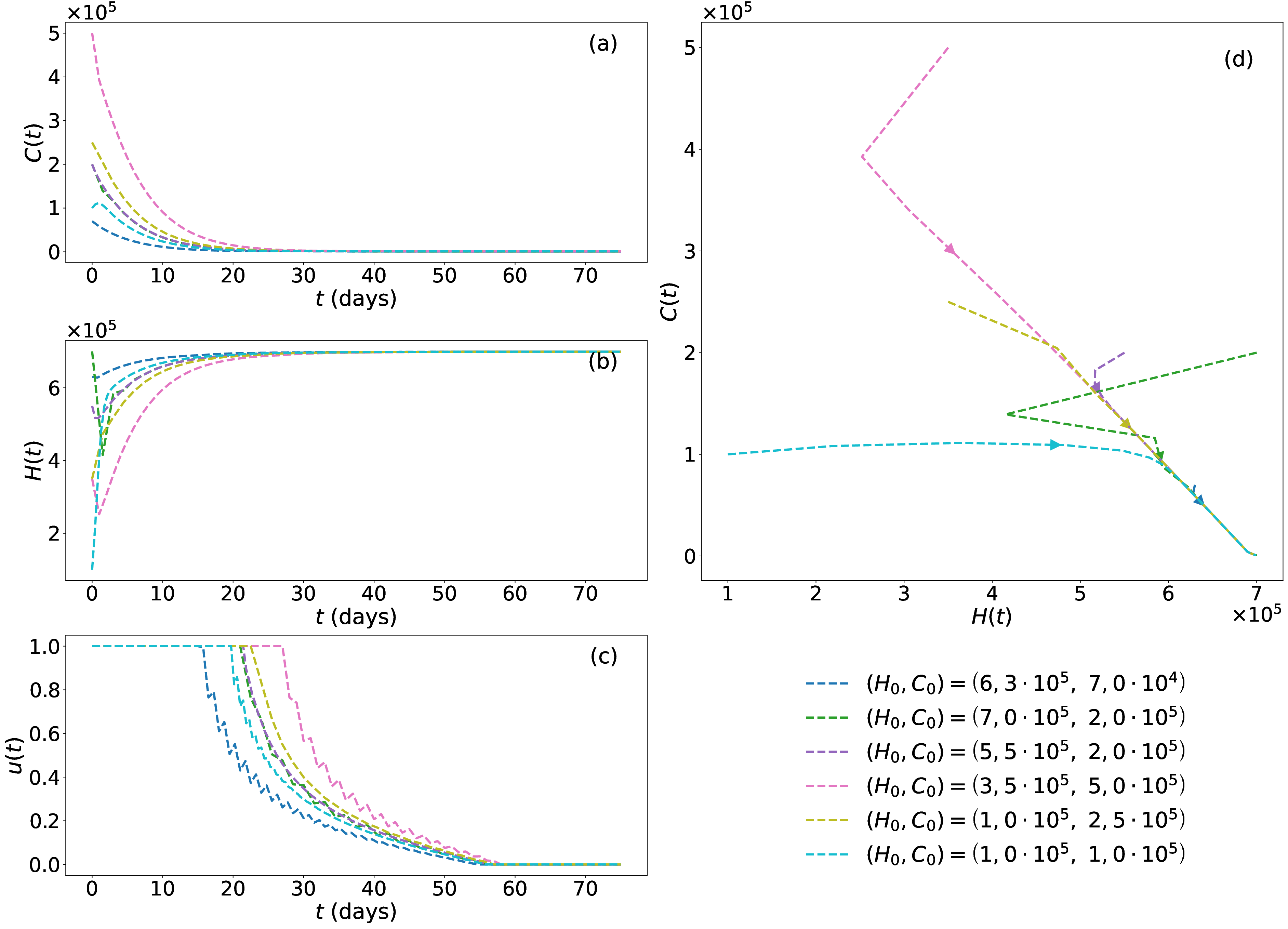}
        \caption{Numerical solution of the optimal control problem \eqref{OCP} for various initial conditions using a direct method. (a) Temporal evolution of the cancer population $C(t)$, showing a rapid and complete elimination within approximately 60 days. (b) Recovery of the healthy population $H(t)$, which effectively returns to the carrying capacity $K = 7 \times 10^5$, bypassing the toxicity-induced displacement observed in the constant dose scenario. (c) Optimal control profiles $u(t)$ exhibiting a ''bang-off" behavior: the treatment maintains maximum intensity $u_{\max}$ initially and then gradually scales down as the tumor is suppressed. (d) Phase portrait $H(t)$ vs. $C(t)$ illustrating the global convergence toward the healthy equilibrium.}
		\label{fig:model_control}
	\end{figure}
	
The numerical results for the optimal control problem are illustrated in Figure \ref{fig:model_control}. A fundamental observation is that the dynamic controller successfully drives the cancer population $C(t)$ to extinction while allowing for a much more robust recovery of the healthy tissue compared to the constant treatment case.

As shown in Figure \ref{fig:model_control}b, the healthy population $H(t)$ steers back toward the carrying capacity $K$ as the treatment $u(t)$ scales down, as shown in Figure \ref{fig:model_control}c. However, it is important to note that within the simulated timeframe, the population does not converge exactly to $K$, but to a high-recovery steady state. This is due to the residual effects of the competitive pressure and the toxicity of the initial high-dose phase. 

The phase portrait in Figure \ref{fig:model_control}d confirms this behavior: all trajectories converge toward the healthy axis ($C=0$), significantly surpassing the shifted equilibrium of the constant dose scenario and demonstrating that the dynamic strategy prioritizes the long-term viability of the host tissue.

\begin{figure}
        \centering
        \includegraphics[width=\textwidth]{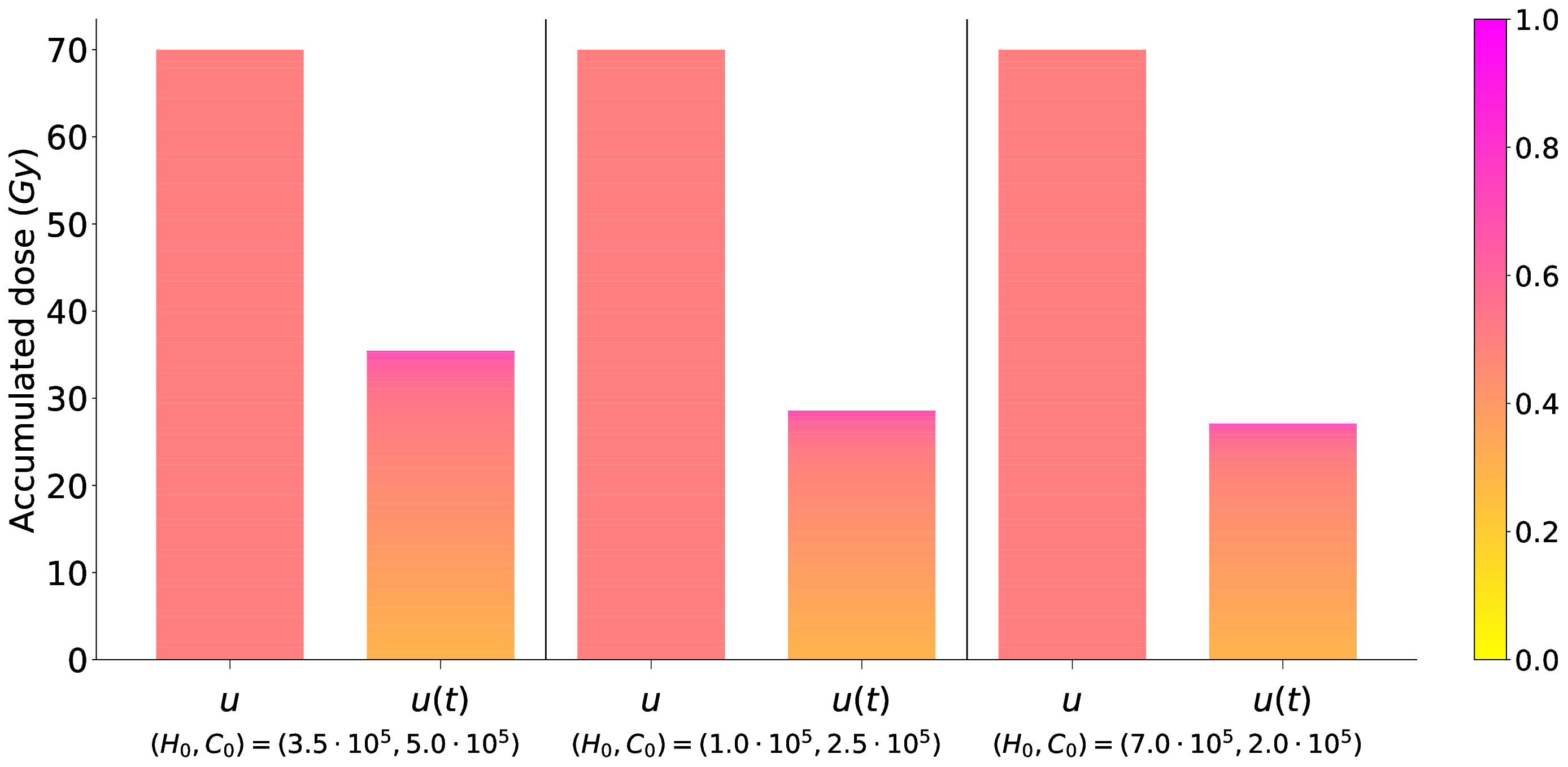}
        \caption{Comparison of total accumulated dose between the constant radiotherapy protocol, $u$, and the optimal control strategy, $u(t)$, for three distinct initial conditions $(H_0, C_0)$. The color intensity, governed by the colorbar, indicates the magnitude of the control variable $u$ at each interval, ranging from minimum, light, to maximum, dark, therapeutic effort.}
        \label{fig:accumulated_dose}
    \end{figure}

The results regarding the total therapeutic effort are summarized in Figure~\ref{fig:accumulated_dose}. A key observation is the difference in dose distribution between the static and dynamic approaches. While the constant treatment $u$ maintains a rigid structure, the optimal strategy $u(t)$ utilizes a time-dependent allocation represented by the segmented bars. The color mapping within these segments, as defined by the associated colorbar, reveals the intensity of the treatment at each stage of the process. Darker tones correspond to periods of maximum radiation intensity, whereas lighter tones indicate a reduction in the dose as the system approaches the desired healthy equilibrium. This adaptive distribution confirms that optimal control not only aims for tumor eradication but also optimizes the timing of dose delivery, potentially reducing the cumulative biological stress on healthy tissue compared to a non-optimized constant protocol.

\section{Conclusions}

The research presented in this manuscript establishes a comprehensive mathematical framework that characterizes the competitive struggle between healthy and cancerous cell populations, moving from basic biological principles to complex therapeutic optimization. By evolving from a single-population logistic growth model to a coupled system of differential equations, we successfully accounted for the metabolic exchanges within the host microenvironment, thereby addressing the limitations of traditional closed-system assumptions. A fundamental contribution of this work was the mathematical demonstration that a shared carrying capacity, in the absence of explicit inter-species competition, results in an unrealistic coexistence. Consequently, the introduction of the competition term $\gamma$ proved essential to reflect the aggressive displacement of healthy tissue observed in clinical reality, confirming that the natural progression of the disease is governed by the Principle of Competitive Exclusion.

The analytical study of the system's stability revealed that the healthy state is inherently unstable when faced with a malignant lineage, meaning that the tumor-dominant equilibrium acts as a terminal attractor for the host. This finding underscores the critical necessity of external control, as the system's autonomous dynamics cannot spontaneously revert to a healthy state. While we demonstrated that constant radiotherapy doses can destabilize the tumor niche, this approach highlights a significant clinical trade-off: the induction of a shifted equilibrium where healthy tissue never fully recovers due to persistent treatment toxicity. 

The application of optimal control theory, implemented through a direct method in \textit{OptimalControl.jl}, represents the most significant achievement of this study. By shifting the paradigm from static treatment to a dynamic, precision-based strategy, we showed that it is possible to navigate the trade-off between tumor eradication and host preservation. Our numerical simulations confirm that a time-dependent, optimized controller is the only approach capable of steering the system's trajectories back to the original carrying capacity $K$, effectively reversing the biological collapse of the host environment. 

Despite these results, this work opens several avenues for future research. The current model could be extended with the inclusion of time delays in the healthy tissue's regenerative response or the consideration of multi-dosage protocols with rest periods could provide even greater clinical depth. Ultimately, this study concludes that integrating dynamic optimal control is not merely a refinement of existing protocols, but a decisive requirement to break the dominance of cancerous lineages and restore biological integrity.

\section*{Declaration of competing interest}
The authors declare that they have no known competing financial interests or personal relationships that could have appeared to influence the work reported in this paper.
	
\section*{Acknowledgments}
\sloppy
Authors gratefully acknowledge financial support by the Spanish Ministerio de Economía y Competitividad and European Regional Development Fund under contract PID 2020-113881RB-I00 AEI/FEDER, UE, and by Xunta de Galicia under Research Grant No. 2021-PG036. All these programs are co-funded by FEDER (UE). The simulations were run in the Supercomputer Center of Galicia (CESGA) and we acknowledge their support.

\bibliographystyle{unsrt}
\bibliography{sample.bib}
	
\end{document}